\documentclass[a4paper,11pt]{article}
\usepackage{paralist}   

\usepackage{amsmath}
\usepackage{amssymb}
\usepackage{amsfonts}
\usepackage[usenames,dvipsnames]{color}

\usepackage{graphicx}
\usepackage{subfig}
\usepackage{bm}
\usepackage{tikz}

\newtheorem{theorem}{Theorem}
\newtheorem{definition}{Definition}
\newtheorem{proposition}[theorem]{Proposition}
\newtheorem{corollary}[theorem]{Corollary}
\newtheorem{lemma}[theorem]{Lemma}
\newtheorem{remark}[theorem]{Remark}
\DeclareMathOperator{\esssup}{ess\,sup}

\newcommand{\ran}{{ Ran}\,}
\renewcommand{\ker}{{ Ker}\,}
\newcommand{\dom}{{ Dom}\,}

\begin{document}

\title{\textbf{Reproducing
pairs and the continuous nonstationary Gabor transform on LCA groups}}
\author{Michael Speckbacher and Peter Balazs\footnote{Acoustics Research Institute, Austrian Academy of Sciences,  Wohllebengasse 12-14, 
1040 Vienna, Austria, speckbacher@kfs.oeaw.ac.at, peter.balazs@oeaw.ac.at}}
\date{}
\maketitle

\vspace{1cm}
\begin{abstract}
In this paper we introduce and investigate the concept of reproducing pairs as a generalization of continuous frames. 
Reproducing pairs yield a bounded analysis and synthesis process while the frame condition can be omitted at both stages.
Moreover, we will investigate certain continuous frames (resp. reproducing pairs) on LCA groups,
which can be described as a continuous version of nonstationary Gabor systems and state sufficient conditions for these systems to 
form a continuous frame (resp. reproducing pair).
As a byproduct we identify the structure of the frame operator (resp. resolution operator).
We will apply our results to systems generated by a unitary action of a subset of 
the affine Weyl-Heisenberg group in $L^2(\mathbb{R})$. This setup will also serve as a nontrivial example of a system for which, whereas
continuous frames exist, no dual system with the same structure exists even if we drop the frame property. 
\end{abstract}
\vspace{2cm}
\textbf{Math Subject Classification:} 22B99, 43A32, 42C15, 42C40.
\textbf{Submitted to:} Journal of Physics A: Mathematical and Theoretical

\newpage

\section{Introduction}
Motivated by physical applications \cite{alanga93a,ka90}, 
Ali et al. \cite{alanga93} and Kaiser \cite{ka94} introduced continuous frames independently in the early 1990's  in order to generalize 
the coherent states approach.
Coherent states are widely used in many areas of theoretical physics, in particular in quantum mechanics, see for instance \cite{alanga00}.
Classically, coherent states are generated by a group action on a single mother wavelet and lead to a resolution of the identity whereas
continuous frames yield a resolution of a positive, bounded, and invertible operator.

Continuous frames have received a lot of interest in recent years \cite{xxlbayasg11,fora05,ranade06},  due to their usefulness in 
signal/image processing \cite{antoine1,duvmus93} as a common background for particular continuous transforms like the wavelet \cite{wei13} 
and other transforms \cite{cado05}. 

There are, however, situations where it is impossible to meet the frame conditions and thus, the goal is to find generalizations 
of this sometimes overly restrictive regime.
To this end one approach is to look at systems that do not satisfy both frame bounds simultaneously leading to the concept of semi-frames 
\cite{jpaxxl11,jpaxxl09,jpaxxl12}.\\

As a next generalizing step, reproducing pairs are introduced in this paper, an approach where a pair of mappings (not necessarily frames,
not even Bessel) in place of a single mapping is used to construct
an invertible analysis/synthesis process.
In the discrete setting the question of dual systems is a current topic of research although restricted to dual frames
\cite{cakula04} or even to dual frames with a particular structure (see  \cite{persampta13,Christensen2014198}). 
The concept of 'weak pairs of dual systems' is for example investigated in  \cite{Weiss1997,felu06,fezi98,JIAHuiFang:1005}. 
In contrast to previous work,  continuous systems are investigated in this paper. Moreover, it will be shown that the property 
of a reproducing pair is not sufficient for the Bessel property. Another difference is that we deal with general Hilbert spaces, not particular function or 
distribution spaces, and, in the first part, systems without a particular structure. 

Reproducing pairs are closely related to frame multipliers \cite{xxlmult1,xxlbayasg11}. In physics, frame multipliers are 
better known as quantization operators \cite{alanga00} and form
the link between classical and quantum mechanics.
Invertibility is a central topic \cite{stoxxl15,balsto09new,uncconv2011} in the mathematical investigation of
multipliers. This includes the question under which conditions a system of two mappings forms a reproducing pair.\\

Reproducing pairs (resp. frames) whose resolution (resp. frame) operator is given by a multiple of the identity are 
of particular interest as inversion is trivial. That raises the following questions: Consider a particular structure and suppose that 
the set of all frames generated by this structure is non-empty. Is there a tight frame contained in this set? If there is no such frame,
can we
find a reproducing pair such that the resolution operator equals the identity? In Section \ref{aWH-section} we show that the answer is in 
general not affirmative.\\

The most common continuous transforms in signal analysis are the short-time Fourier transform (STFT) \cite{groe1} and the continuous wavelet transform (CWT)
\cite{ma09}. In applications they are used in their sampled, discretized version. 
Both transforms have a time-frequency resolution which is either fixed for all frequencies (for the STFT) or follows a given rule 
(for wavelets). 
In practice, functions or signals often show particular time-frequency characteristics which call for adaptive and adaptable representations 
\cite{badokowto13}. In \cite{nsdgt10} the authors introduced adaptivity either in time or frequency where perfect reconstruction is still
possible. We will modify this approach within this paper and introduce a continuous version of nonstationary Gabor frames which 
is conceptually a version of generalized shift-invariant systems \cite{helawe02,rosh05}.

To achieve this, we will regard systems on $L^2(X,\mu)$ where we intend to choose the domain $X$ such that on the one hand, $X$ is
as general as
possible and, on the other hand, $X$ possesses enough structure such that the definition of nonstationary Gabor systems
is meaningful and fruitful.
It turns out that locally compact abelian (LCA) groups are the appropriate domains for our investigation since Fourier analysis on 
LCA groups diagonalizes convolution in the right function spaces.
We will provide a substitute for the frame (resp. reproducing pair) condition and show that the frame (resp. resolution) operator is
given by a Fourier multiplier.

The study of representations of the affine Weyl-Heisenberg  group $G_{aWH}$ is of particular interest since 
both, the Weyl-Heisenberg group and the affine group which are the underlying groups of the STFT and the CWT respectively, are subgroups 
of $G_{aWH}$. Transforms generated by an appropriate restriction of $G_{aWH}$ form a subclass
 of continuous nonstationary Gabor transforms and we will apply the results to this setting. In particular, the $\alpha$-transform 
 (see \cite{daforastte08,dalo06}) a class 
 of transforms intermediate between 
 the STFT and the CWT will serve as an example in this paper. \\

The paper at hand is organized as follows: in Section \ref{prel} we will briefly present the basic results on Fourier analysis on LCA
groups and continuous frames.  In Section \ref{sec:reporpair0} we introduce 
the concept of reproducing pairs. Section \ref{cont-nonstationary gabor} is concerned with the continuous nonstationary Gabor 
transform on LCA groups. The results are then applied in Section \ref{aWH-section}
to representations of subsets of the affine Weyl-Heisenberg group and we will show particularities of reproducing 
pairs in comparison to continuous frames with the help of an example.

\section{Preliminaries}\label{prel}

Throughout this paper we will assume that any mapping 
$\Psi:X\rightarrow \mathcal{H}$ indexed by some measure space $(X,\mu)$ is essentially norm-bounded, 
i.e., $\esssup_{x\in X}\|\Psi(x)\|_\mathcal{H}\leq C$ and  the space of bounded linear operators
 with bounded inverse from $\mathcal{H}$ to $\mathcal{K}$ will be denoted by  $GL(\mathcal{H},\mathcal{K})$.

\subsection{Fourier analysis on LCA groups}

Fourier analysis is the most important tool in harmonic analysis and has been extended in mid-twentieth century to functions on
locally compact topological groups. In this paper we will in particular use Fourier theory on LCA groups.
For a thorough introduction consult the standard text books \cite{fo95,lo53}. The  fundamental examples of 
LCA groups in harmonic analysis are the additive groups $\mathbb{R},\ \mathbb{Z},\ \mathbb{R}/\mathbb{Z}$ and $\mathbb{Z}/ N\mathbb{Z}$ and 
their d-fold products. Their relation is depicted in the following diagram, see e.g. \cite{so05}.
\begin{center}
\begin{tikzpicture}[scale=0.7]
\draw[->,thick] (0,1.8)--(0,1)node[pos=0.5,left,font=\footnotesize]{Sampling};
\draw[->,thick] (0.5,0.5)--(3.3,0.5)node[pos=0.45,below,font=\footnotesize]{Periodization};
\draw[->,thick] (0.5,2.3)--(3.3,2.3)node[pos=0.45,above,font=\footnotesize]{Periodization};
\draw[->,thick] (3.9,1.8)--(3.9,1)node[pos=0.5,right,font=\footnotesize]{Sampling};
\draw (0,0.5) node {$\mathbb{Z}$};
\draw (0,2.3) node {$\mathbb{R}$};
\draw (3.9,2.3) node {$\mathbb{R}/\mathbb{Z}$};
\draw (4.1,0.5) node {$\mathbb{Z}/N\mathbb{Z}$};
\end{tikzpicture}
\end{center}
Every LCA group possesses a unique translation invariant 
measure on $G$ (up to a constant factor) called the Haar measure  $dx$. Convolution of two functions is given by
$f\ast g(y):=\int_G f(x)g(x^{-1}y)dx$.
It follows by the Riesz-Thorin theorem that convolution with a fixed function $g\in L^1(G)$
is a bounded operator on $L^p(G)$ as $\|f\ast g\|_p\leq \|f\|_p\|g\|_1$ for $1\leq p\leq\infty$.

A character $\xi$ is a continuous homomorphism from $G$ to the torus $\mathbb{T}$, i.e., $\xi(xy)=\xi(x)\xi(y)$ and $|\xi(x)|=1$. 
The set $\widehat G$ of all characters of $G$ is called the dual group and is again an LCA group with pointwise multiplication and the 
topology of compact convergence on $G$. 
The Pontryagin duality theorem states that any LCA group is reflexive, i.e., the dual group of $\widehat G$ is isomorphic 
to $G$. The dual groups of the fundamental examples are given by
$\widehat{\mathbb{R}}\cong\mathbb{R},\  \widehat{\mathbb{R}/\mathbb{Z}}\cong\mathbb{Z},\ 
 \widehat{\mathbb{Z}}\cong\mathbb{R}/\mathbb{Z}$ and  $\widehat{\mathbb{Z}/N\mathbb{Z}}\cong\mathbb{Z}/N\mathbb{Z}$.
 
Now we are able to define the Fourier transform on $L^1(G)$ by
\begin{equation*}
 \hat f(\xi):=\int_G \overline{\xi(x)}f(x)dx,\ \ \xi\in\widehat G
\end{equation*}
It can be shown that this definition extends to an isometric isomorphism from $L^2(G)$ to $L^2(\widehat G)$ if the Haar
measure on $\widehat G$ is appropriately normalized, i.e., 
$\|f\|_2=\|\hat f\|_2\ \forall f\in L^2(G)$, and that Parseval's formula $\langle f,g\rangle=\langle \hat f,\hat g\rangle$ holds for all
$f,g\in L^2(G)$. 
In addition, if $f,g\in L^2(G)$ and $f\ast g\in L^2(G)$ it follows that $(f\ast g)\ \widehat{}\ (\xi)=\hat f(\xi) \hat g(\xi)$.

\subsection{Frames} 
\noindent Frames were first introduced by Duffin and Schaeffer \cite{duscha52} as a  generalization of orthonormal bases. 
Ali et al.  \cite{alanga93} and Kaiser \cite{ka90} extended  frames to mappings on a measure space $(X,\mu)$.
\begin{definition}\label{def-cont-frame}
Let $\mathcal{H}$ be a Hilbert space and $(X,\mu)$ be a measure space.
A mapping $\Psi:X\rightarrow \mathcal{H}$  is called a continuous frame if
\begin{enumerate}[(i)]
 \item $\Psi$ is weakly measurable, i.e., $x\mapsto\langle f,\Psi(x)\rangle$ is a measurable function for all $f\in\mathcal{H}$ 
 \item there exist positive constants $A,B>0$ s.t. 
 \begin{equation}\label{frame-condition}
  A\left\|f\right\|_\mathcal{H}^2\leq\int_{X}\left|\langle f,\Psi(x)\rangle\right|^2d\mu(x)\leq B\left\|f\right\|_\mathcal{H}^2,\ \
  \forall f\in\mathcal{H}
 \end{equation}
\end{enumerate}
\end{definition}
The mapping $\Psi$ is called Bessel if the right inequality in (\ref{frame-condition}) is satisfied.
The definition of  a continuous frame concurs with the standard definition of a  frame if   $X$  is a countable set and $\mu$ is the counting
measure, see for example \cite{christ1}. For a short and self-contained introduction to continuous frames, see \cite{ranade06}.

Let us define the basic operators in frame theory: the analysis operator
\begin{equation*}
 C_\Psi:\mathcal{H}\rightarrow L^2(X,\mu),\ \ \ C_\Psi f(x):=\langle f,\Psi(x)\rangle
\end{equation*}
and the synthesis operator
\begin{equation*}
 D_\Psi:L^2(X,\mu)\rightarrow \mathcal{H},\ \ \ D_\Psi F:=\int_X F (x)\Psi(x)d\mu(x)
\end{equation*}
where the integral is defined weakly. Observe that $C_\Psi^\ast =D_\Psi$ if $\Psi$ is Bessel. The frame operator is defined by composition
of $C_\Psi$ and $D_\Psi$
\begin{equation*}
 S_\Psi:\mathcal{H}\rightarrow \mathcal{H},\ \ \ S_\Psi f:=D_\Psi C_\Psi f=\int_X\langle f,\Psi(x)\rangle \Psi(x)d\mu(x)
\end{equation*}
$S_\Psi$ is self-adjoint, positive, bounded and invertible. The canonical dual $S^{-1}_\Psi\Psi$ also forms a frame
with frame bounds $B^{-1},A^{-1}$, and the following reproducing formula holding weakly
\begin{equation}
f=D_{S^{-1}_\Psi\Psi}C_\Psi f=D_\Psi C_{S^{-1}_\Psi\Psi}f,\ \ \ \forall f\in\mathcal{H}
\end{equation}

The analysis operator $C_\Psi$ is in general not onto $L^2(X,\mu)$ since a continuous frame $\Psi$  is overcomplete
whenever the underlying measure space is not atomic \cite{hedera00}.
The range of the $C_\Psi$ is nevertheless characterized by the reproducing kernel of $\Psi$, i.e., for $F\in L^2(X,\mu)$
there exists $f\in\mathcal{H}$ s.t. $F=C_\Psi f$ if and only if
$F(x)=R(F)(x)$ where $R$ is an integral operator with kernel $\mathcal{R}(x,y):=\langle S_\Psi^{-1} \Psi(y),\Psi(x)\rangle$ and
\begin{equation*}
 R(F)(x):=\int_X\mathcal{R}(x,y)F(y)d\mu(y)
\end{equation*}
In addition, $R$ is the orthogonal projection from $L^2(X,\mu)$ onto $\ran(C_\Psi)$.

Since, in general, $\ker(D_\Psi)\neq\{0\}$, there may exist other mappings $\Psi^d$ different from $S_\Psi^{-1}\Psi$ satisfying
\begin{equation*}
f=D_{\Psi^d}C_\Psi f=D_\Psi C_{\Psi^d}f,\ \ \ \forall f\in\mathcal{H}
\end{equation*}
Such a mapping $\Psi^d$ is called a dual frame.

\section{Reproducing Pairs} \label{sec:reporpair0}
Dual frames have been studied in various articles, see for example 
\cite{cakula04} or  \cite{Christensen2014198,persampta13,frgawawe97} for Gabor and wavelet frames. Nevertheless, it has been pointed 
out in \cite{jpaxxl09} that it is sometimes not possible to satisfy both frame conditions at the same time. Therefore, the authors 
introduced upper/lower
semi-frames, i.e., total systems that only satisfy the upper/lower frame inequality. 
In  the paper at hand, we will present a different approach  as a generalization of the concept of weakly dual pair of systems 
\cite{Weiss1997,felu06} motivated by 
\cite[Definition 1.1.1.33]{te01} which combines expansion via two reproducing mappings and the omission of the frame property.
It is clear that even for the discrete setting dual systems exist which are non-Bessel, see e.g.  \cite{bstable09} Example 4.1.7(i).

\begin{definition}\label{rep-pair-definition}
 Let $(X,\mu)$ be a measure space and $\Psi,\Phi:X\rightarrow\mathcal{H}$ weakly measurable.
 The pair of mappings $(\Psi,\Phi)$ is called a reproducing pair for $\mathcal{H}$ if the resolution operator 
 $S_{\Psi,\Phi}:\mathcal{H}\rightarrow \mathcal{H}$ weakly defined by
 \begin{equation}\label{rep-pair-def}
 S_{\Psi,\Phi} f:=\int_X \langle f,\Psi(x)\rangle \Phi(x)d\mu(x)
 \end{equation}
is an element of $GL(\mathcal{H})$.
\end{definition}
 Note that if $(\Psi,\Phi)$ is a reproducing pair, $(\Psi,  S_{\Psi,\Phi}^{-1} \Phi)$ is a dual system. So the results in here can be 
 reinterpreted as results about dual systems.

This definition is indeed a generalization of continuous frames because on the one hand neither $\Psi$ nor $\Phi$ are required 
to meet the frame condition.
On the other hand, a weakly measurable mapping $\Psi$ is a continuous frame if and only if $(\Psi,\Psi)$ is a reproducing pair. This 
feature follows immediately if one considers the original definition of  a continuous frame in \cite{alanga93} which is equivalent to 
Definition \ref{def-cont-frame}.

In the very same article the authors showed that it is possible to generate new, equivalent continuous frames  from a given continuous
frame. Within the present paper, we will use a more general concept of equivalent mappings from \cite{jale14} and adapt it to reproducing pairs. 

 \begin{lemma}\label{reprod-equiv-lemma}
Let $\mathcal{H},\ \mathcal{K}$ be  Hilbert spaces, $(X,\mu),\ (Y,\nu)$ be a measure space, $\rho:Y\rightarrow X$ a bijective mapping that
satisfies $\nu\circ \rho^{-1}=\mu$ and preserves measurability,
 $T\in GL(\mathcal{H},\mathcal{K})$ and $\tau:Y\rightarrow \mathbb{C}$ a measurable function with $|\tau(y)|=1$.\\
  Define $\widetilde \Psi(y):=\tau(y) T ( \Psi\circ \rho) (y)$ (and $\widetilde \Phi$ respectively), 
  then $(\Psi,\Phi)$ is a reproducing pair for $\mathcal{H}$ with respect to $(X,\mu)$ if and only if $(\widetilde \Psi,\widetilde \Phi)$ 
  is a reproducing pair for $\mathcal{K}$ with respect to $(Y,\nu)$.
 \end{lemma}
 \textbf{Proof:} Let $f,g\in\mathcal{K}$. It holds
 \begin{align*}
  \langle S_{\widetilde \Psi,\widetilde \Phi} f,g\rangle_{\mathcal{K}} &=\int_{Y}\big\langle f,\tau(y) T ( \Psi\circ \rho) (y)\big
  \rangle_{\mathcal{K}}\big \langle \tau(y) T ( \Phi\circ \rho) (y),g\big\rangle_{\mathcal{K}} d\nu(y)\\
&=\int_{Y}\big\langle T^\ast f, (\Psi\circ \rho) (y)\big\rangle_{\mathcal{H}}\big\langle (\Phi\circ \rho) (y),T^\ast g\big\rangle_{\mathcal{H}} d\nu(y)\\
&=\int_{X}\big\langle T^\ast f,\Psi(x)\big\rangle_{\mathcal{H}}\big\langle \Phi(x),T^\ast g\big\rangle_{\mathcal{H}} d\mu(x)\\
&=\langle S_{ \Psi,\Phi} T^\ast f,T^\ast g\rangle_{\mathcal{H}}\\
&=\langle TS_{ \Psi,\Phi} T^\ast f, g\rangle_{\mathcal{K}}
\end{align*}
 Hence, we can identify the resolution operator $S_{\widetilde \Psi,\widetilde \Phi} = TS_{ \Psi,\Phi} T^\ast$ and the result follows as 
 $S_{\widetilde \Psi,\widetilde \Phi}\in GL(\mathcal{K})$ if and only if
 $S_{ \Psi,\Phi}\in GL(\mathcal{H})$.\hfill $\Box$

\vspace{0.3cm}
\noindent Unlike the frame operator $S_\Psi$, $S_{\Psi,\Phi}$ is in general neither positive nor self-adjoint, since $S_{\Psi,\Phi}^\ast=
S_{\Phi,\Psi}$.\\ In the following, we need to define the domain of $D_\Phi$ as
\begin{equation*}
\dom(D_\Phi):=\Big\{F\in L^\infty(X,\mu):\ \int_X F(x)\Phi(x)d\mu(x)\ converges\ weakly\Big\}
\end{equation*}
It is important to notice here that this definition differs from the standard definition, where  $\dom(D_\Phi)$ is considered to be
 a subspace of $L^2(X,\mu)$. This is justified since $\Psi$ is not necessarily Bessel and consequently $\ran(C_\Psi)\nsubseteq L^2(X,\mu)$ 
 in general but $\ran(C_\Psi)\subseteq \dom(D_\Phi)\subseteq L^\infty(X,\mu)$ by the standing assumption on $\Psi$.
Using the standard definition of $\dom(D_\Phi)$ leads to a possibly unbounded resolution operator and new interesting perspectives 
for future research on reproducing pairs.

We will now give a necessary condition for $F:X\rightarrow \mathbb{C}$ to be an element of the 
image of $C_\Psi$ in terms of a reproducing kernel.

\begin{proposition} \label{sec:repordkern1}
 Let $(\Psi,\Phi)$ be a reproducing pair for $\mathcal{H}$ and $F\in \dom(D_\Phi)$. It holds that $F(x)=\langle f,\Psi(x)\rangle$, 
 for almost every $x\in X$ and some $f\in\mathcal{H}$, if and only if $F(x)=R(F)(x)$ with the integral kernel
 \begin{center}
 $\mathcal{R}(x,y)=\langle S_{\Psi,\Phi}^{-1}\Phi(y),\Psi(x)\rangle$
 \end{center}
 Moreover, $L^1(X,\mu)\cap L^\infty(X,\mu)\subset \dom(D_\Phi)$, which in particular implies that $\dom(D_\Phi)\cap L^2(X,\mu)$ is 
 dense in $L^2(X,\mu)$.
\end{proposition}
\textbf{Proof:}
Let $F(x)=\langle f,\Psi(x)\rangle$, then 
\begin{align*}
 R(F)(x) &=\int_X \langle f,\Psi(y)\rangle\langle \Phi(y),(S_{\Psi,\Phi}^{-1})^\ast\Psi(x)\rangle d\mu(y)\\
&=\langle S_{\Psi,\Phi}f, (S_{\Psi,\Phi}^{-1})^\ast\Psi(x)\rangle=C_\Psi f(x)=F(x)
\end{align*}
Assume now that $R(F)(x)=F(x)$. Since $F\in \dom(D_\Phi)$
we can weakly define 
$g:=\int_{X}F(x)\Phi(x)d\mu(x)\in\mathcal{H}$. It then follows that $F(x)=C_\Psi f(x)$ where $f:=S_{\Psi,\Phi}^{-1}g$, since
\begin{align*}
 C_\Psi f(x) &=\langle S_{\Psi,\Phi}^{-1}g,\Psi(x)\rangle=\langle g,(S_{\Psi,\Phi}^{-1})^\ast\Psi(x)\rangle\\
 &=\int_XF(y)\langle \Phi(y),(S_{\Psi,\Phi}^{-1})^\ast\Psi(x)\rangle d\mu(y)=R(F)(x)=F(x)
\end{align*}
It remains to show that if $F\in L^1(X,\mu)\cap L^\infty(X,\mu)$ it follows that the integral $\int_{X}F(x)\Phi(x)d\mu(x)$ converges weakly.
Let $h\in\mathcal{H}$
\begin{align*}
\big|\big\langle \int_{X}F(x)\Phi(x)d\mu(x),h\big\rangle\big| &\leq\int_{X}\big|F(x)\langle\Phi(x),h\rangle \big|d\mu(x)\\
 &\leq \sup\limits_{x\in X}\|\Phi(x)\|_\mathcal{H}\ \|F\|_{L^1(X,\mu)}\ \|h\|_\mathcal{H}
\end{align*}
and hence by Riesz representation theorem $F\in \dom(D_\Phi)$.
\hfill $\Box$

\vspace{0.3cm}
Notice that, unlike for frames, there may exist $f\in\mathcal{H}$ s.t.
$C_\Psi f \notin L^2(X,\mu)$. An example for such a setting can be found in Section \ref{example-gawh-section}.
However, it is not difficult to see that the images of $C_\Psi$ and $C_\Phi$ are subspaces of mutually dual spaces with the duality pairing 
$\langle F,H\rangle=\int_X F(x)\overline{H}(x)d\mu(x)$.\\

In the course of investigating reproducing pairs and continuous frames generated by a particular structure, for example by the action of 
a group representation  on a single window, three questions naturally arise:
\begin{enumerate}[(i)]
 \item Are there equivalent or sufficient conditions for $\Psi$ (resp. $(\Psi,\Phi)$) to be a 
continuous frame (resp. a  reproducing pair)? 
\item What can be said about the structure of the frame operator (resp. resolution operator)?
\item Given a ``nice'' structure which generates $\Psi$ and $\Phi$: Can one choose  $\Psi$ and $\Phi$ such that 
the resolution operator is the identity operator?
\end{enumerate}

In \cite{grmopa86} the authors gave a sufficient answer to these questions if $X=G$ is a locally compact group and
$\mu$ its left Haar measure. If $\pi:G\rightarrow\mathcal{H}$ is a square-integrable group representation, i.e.,  if it is irreducible and
\begin{equation*}
\mathcal{A}:=\Big\{\psi\in\mathcal{H}:\ \int_G\left|\langle \pi(x)\psi,\psi\rangle\right|^2d\mu(x)<\infty\Big\}\neq\{0\},
\end{equation*}
then there exists a unique self-adjoint operator $L$ with dense domain $\mathcal{A}$ s.t. for all
$\psi,\varphi\in\mathcal{A}$ the following orthogonality relation holds
\begin{equation*}
 \int_G\langle f_1,\pi(x)\psi\rangle\overline{\langle f_2,\pi(x)\varphi\rangle}d\mu(x)=\langle L\varphi, L\psi\rangle\langle f_1,f_2\rangle,
 \ \ \ \forall f_1,f_2\in \mathcal{H}
\end{equation*}
Elements of $\mathcal{A}$ are called admissible windows. If $G$ is unimodular, then $L$ is a multiple of the identity. Hence, we see that,
after normalization,
the resolution operator is given by the identity operator. 

Nevertheless, the restriction  to systems arising from square-integrable group representations excludes a wide range of interesting
transforms. This is why we will introduce more flexible transforms in the following.

\section{The continuous nonstationary Gabor transform on LCA groups}\label{cont-nonstationary gabor}
In this section we consider continuous systems on $L^2(G)$ motivated by nonstationary Gabor frames. 
Such systems were first introduced as generalized shift-invariant systems, see for example \cite{helawe02,rosh05,christel05} and the parallel work \cite{jale14},
in order to gain flexibility in analyzing signals with specific time-frequency characteristics.
The main focus of these papers,  however, is to find sufficient conditions for tight frames and mutually dual frames.
 Balazs et al. \cite{nsdgt10} introduced the terminology of nonstationary Gabor systems, motivated from a signal-processing viewpoint. 
 The very same study concentrates on the frame property and 
easy inversion via painless nonorthogonal expansions (see \cite{daubgromay86}), i.e., frame expansions whose frame operator is diagonal or 
diagonal in the Fourier domain.

\subsection{Translation invariant systems}

Throughout the rest of this paper we assume $\mathcal{H}=L^2(G)$, where $G$ is a second countable LCA group. In particular, this assumption
implies that $L^2(G)$ is separable, both $G$ and $\widehat G$ are $\sigma$-compact and  the Haar measures 
on $G$ and $\widehat G$ are, consequently, $\sigma$-finite.
The translation operator on $G$ is given by $T_z f(x):=f(z^{-1}x),\ x,z\in G$ and its Fourier transform 
is given by $\widehat{T_z f}(\xi)=\xi(z^{-1})\hat f(\xi),\ \xi\in\widehat G$.
Now let $\psi_y, \varphi_y\in L^2(G)$ for all $y\in Y$ with $(Y,\mu)$ being a $\sigma$-finite measure space. For  
$(x,y)\in G\times Y$ we define
\begin{equation*}
 \Psi(x,y):=T_{x}\psi_y\hspace{0.5cm} and \hspace{0.5cm}  \Phi(x,y):=T_{x}\varphi_y
 \end{equation*}
 and the continuous nonstationary Gabor transform (CNSGT) by
 \begin{equation*}
 C_\Psi f(x,y):=\langle f,\Psi(x,y)\rangle 
 \end{equation*}
{The following theorem gives sufficient condition for $(\Psi,\Phi)$ to be a reproducing pair and reveals the structure of the reproducing (resp. frame) operator.
A related result, where $\Psi,\Phi$ are assumed to be Bessel, can be found in the parallel work \cite{jale14}.

\begin{theorem}\label{cnsgt}
If there exist $A, B, C>0$, s.t.
 \begin{equation}\label{rep-pair-con}
  A\leq |m_{\Psi,\Phi}(\xi)|\leq B,\ \text{for a.e.}\ \xi\in\widehat{G}
 \end{equation} 
with
  \begin{equation}
m_{\Psi,\Phi}(\xi):=\int_Y \overline{\widehat{\psi_y}(\xi)}\widehat{\varphi_y}(\xi)d\mu(y)
 \end{equation}
 and
   \begin{equation}\label{L1-con}
\int_Y \big|\widehat{\psi_y}(\xi)\widehat{\varphi_y}(\xi)\big|d\mu(y)\leq C,\ \text{for a.e.}\ \xi\in\widehat{G}
 \end{equation}
 then $(\Psi,\Phi)$ is a 
 reproducing pair for $L^2(G)$ and resolution operator is given weakly by
 \begin{equation}\label{res-op-rep}
  S_{\Psi,\Phi} f=\mathcal{F}^{-1}(m_{\Psi,\Phi}\cdot\mathcal{F}(f))
 \end{equation}
 If $\Psi=\Phi$, then $\Psi$ is Bessel if and only if the upper bound in (\ref{rep-pair-con}) is satisfied and a continuous frame with 
 frame operator $S_{\Psi}=S_{\Psi,\Psi}$ and frame bounds $A,B$ if and only if condition (\ref{rep-pair-con}) is satisfied. In particular, 
 the frame is tight if $A=B$.
\end{theorem}
\textbf{Proof:} Let $f_1,f_2\in 
L^1(G)\cap L^2(G)$, $\psi_y,\varphi_y\in L^2(G)$ and assume that (\ref{rep-pair-con}) and (\ref{L1-con}) hold.
Observe that $\langle f,T_{x}\psi_y\rangle=f\ast \psi_y^\ast(x)$, where $g^\ast(x):=\overline{g}(x^{-1})$ is the involution of $g$. 
Since $f\in L^1(G)$ it follows that $f\ast \psi_y^\ast\in L^2(G)$ and therefore $(f\ast \psi_y^\ast)\ \widehat{}\ (\xi)=\hat f(\xi)
\overline{\widehat \psi_y(\xi)}$.
Parseval's formula yields
\begin{align*}
 \langle S_{\Psi,\Phi} f_1,f_2\rangle&=\int_Y\int_G\langle f_1,T_{x}\psi_y\rangle\overline{\langle f_2,T_{x}
 \varphi_y\rangle} dxd\mu(y)\\
&=\int_Y\int_{\widehat{G}}\hat{f}_1(\xi)\overline{\hat{f}}_2(\xi)
\overline{\widehat{\psi_y}(\xi)}\widehat{\varphi_y}(\xi)d\xi d\mu(y)\\
 &=\int_{\widehat G}m_{\Psi,\Phi}(\xi)\hat{f}_1(\xi)\overline{\hat{f}}_2(\xi)d\xi\\ &=\langle \mathcal{F}^{-1}(m_{\Psi,\Phi}\cdot\mathcal{F}
 (f_1)),f_2\rangle
\end{align*}
where condition (\ref{L1-con}) guarantees that Fubini's theorem is applicable. 
$S_{\Psi,\Phi}$  extends to a bounded and invertible operator on $L^2(G)$ by a standard density argument and (\ref{rep-pair-con}).

It remains to show that if $\Psi$ is Bessel, it follows that the frame operator is given by (\ref{res-op-rep}). 
By the previous calculation we get
\begin{equation*}
 \langle S_{\Psi} f,f\rangle=\int_Y\int_{\widehat G}|\hat f(\xi)|^2
|\widehat{\psi_y}(\xi)|^2 d\xi d\mu(y)\leq B\|f\|^2_2
\end{equation*}
Consequently, Fubini's theorem is again applicable and the frame operator is given by $S_\Psi f=\mathcal{F}^{-1}(m_{\Psi}\cdot
\mathcal{F}(f))$. It is easy to see that $S_\Psi$ is bounded with bounded inverse only if the symbol
$m_\Psi$ is essentially bounded from above and below.\hfill $\Box$

\vspace{0.3cm}
With a slight misuse of terminology we call $\{\psi_y\}_{y\in Y}$ admissible if (\ref{rep-pair-con}) is satisfied. 
If conditions 
(\ref{rep-pair-con}) and (\ref{L1-con}) are satisfied,
 $\{(\psi_y,\varphi_y)\}_{y\in Y}$ is called  cross-admissible.
Note that the inverse of a Fourier multiplier is given 
by another Fourier multiplier with the inverse symbol, i.e., $S_{\Psi,\Phi}^{-1}f=\mathcal{F}^{-1}(m_{\Psi,\Phi}^{-1}\cdot\mathcal{F}(f))$.

\begin{remark} \rm
  The property that  $m_{\Psi,\Phi}(\xi)\in\mathbb{C}\backslash\{0\}$
  reveals why $S_{\Psi,\Phi}$ is in general neither self-adjoint nor positive whereas  $m_\Psi(\xi)\in \mathbb{R}_{>0}$
  guarantees self-adjointness and positivity of $S_{\Psi}$.
\end{remark}

\begin{corollary}\label{cor-trans-inv}
Suppose that $S_{\Psi,\Phi}$ is given by (\ref{res-op-rep}). $S_{\Psi,\Phi}=Id$ if and only if $m_{\Psi,\Phi}(\xi)=1$, for a.e.
$\xi\in\widehat G$.\\
The canonical dual of a translation invariant frame $\Psi$ is another translation invariant system $\Phi(x,y):=T_xS_\Psi^{-1}\psi_y$
\end{corollary}
\textbf{Proof:} 
Since the Fourier transform is an isometric isomorphism it follows that the equation  $f=\mathcal{F}^{-1}(m_{\Psi,\Phi}^{-1}\cdot\mathcal{F}(f))$
 holds for all $f\in L^2(G)$ if and only if the Fourier transform of $f$ is not altered in $L^2(G)$, i.e., $m_{\Psi,\Phi}(\xi)=1$, for 
 a.e. $\xi\in\widehat G$.
 
 To proof the last assertion we only have to show that the translation operator commutes with the inverse frame operator. This is 
 obviously the case since the inverse frame operator is a Fourier multiplier operator and translation corresponds to character 
 multiplication in Fourier domain.
\hfill $\Box$

\subsection{Character invariant systems}

Instead of shifting $\psi_y$ along $G$ we will now multiply the windows $\psi_y$ with a character $\xi\in\widehat G$, i.e., we consider the 
modulation operator $M_\xi f(x):=\xi(x)f(x)$ and the mappings 
\begin{equation*}
 \Psi(\xi,y):=M_{\xi}\psi_y\hspace{0.5cm}  and \hspace{0.5cm}   \Phi(\xi,y):=M_{\xi}\varphi_y
 \end{equation*}
 where  $(\xi,y)\in\widehat{G}\times Y$ and derive a similar result as in Theorem \ref{cnsgt}.
\begin{corollary} 
The pair of mappings $(\Psi,\Phi)$ is a reproducing pair for  $L^2(G)$, if there exist 
 $A, B, C>0$ s.t.
 \begin{equation}\label{sym-con2}
  A\leq |m_{\Psi,\Phi}(x)|\leq B,\ \text{for a.e.}\ x\in G
 \end{equation} 
 where
 \begin{equation}
 m_{\Psi,\Phi}(x):=\int_Y \overline{\psi_y(x)}\varphi_y(x)d\mu(y)
 \end{equation}
  and
   \begin{equation}
\int_Y \big|\psi_y(x)\varphi_y(x)\big|d\mu(y)\leq C,\ \text{for a.e.}\ x\in G
 \end{equation}
 The resolution operator is weakly given by
 \begin{equation}
  S_{\Psi,\Phi} f=m_{\Psi,\Phi}\cdot f
 \end{equation}
 If $\Psi=\Phi$, then $\Psi$ is Bessel if and only if the upper bound in (\ref{sym-con2}) is satisfied. Moreover, $\Psi$ is
 a continuous frame with frame operator $S_{\Psi}=S_{\Psi,\Psi}$ and frame bounds $A,B$ if and only if condition
 (\ref{sym-con2}) is satisfied. In particular, the frame is tight if $A=B$ in (\ref{sym-con2}).
\end{corollary}
\textbf{Proof:} Using that $\widehat{M_\omega f}(\xi)=T_\omega\hat f(\xi)$  implies 
$\langle f,M_{\xi} \psi_y\rangle=\langle \hat{f},T_{\xi} \widehat{\psi_y}\rangle$. The same steps as in the proof of Theorem \ref{cnsgt} 
together with  $\mathcal{F}_{\widehat{G}}\mathcal{F}_Gf(x)=f(x^{-1})$ conclude the proof.\hfill $\Box$

\subsection{Examples}\label{examples}
Let us apply previous results to two short examples with $G=(\mathbb{R},+)$. Notice that in this situation $\widehat G\cong G$.
\paragraph{Example 4.3. (a)} For the short-time Fourier system
$\Psi(x,\omega)=M_{\omega}T_x\psi$, $\Phi(x,\omega)=M_{\omega}T_x\varphi$, $(x,\omega)\in \mathbb{R}^{2d}$ with $\mu$ the Lebesgue measure, 
one gets the Fourier symbol 
$m_{\Psi,\Phi}(x)=\langle T_x\varphi,T_x\psi\rangle=\langle \varphi,\psi\rangle$ and the well-known inversion 
formula
\begin{equation*}
 f=\frac{1}{\langle \varphi,\psi\rangle}\int_{\mathbb{R}^{2d}}\langle f,M_\omega T_x \psi\rangle M_\omega T_x\varphi dxd\omega
\end{equation*}
\paragraph{Example 4.3. (b)} Here we show that the theory also applies to discrete measure spaces with weighted counting measure. 
Consider the semi-discrete wavelet system on a dyadic scale grid, i.e.,
$\Psi(x,j)=T_xD_{2^j}\psi,\ (x,j)\in \mathbb{R}\times\mathbb{Z}$, with $D_a$ denoting the dilation operator $D_af(x):=a^{-1/2}f(x/a)$, and $Y=\mathbb{Z}$ 
equipped with the weighted counting measure $\mu(j)=2^{-j}$. The Fourier symbol then reads
\begin{equation*}
 m_{\Psi}(\xi)=\sum_{j\in\mathbb{Z}}|\hat\psi(2^j\xi)|^2
\end{equation*}
It is not difficult to verify that  the symbol $m_{\Psi}$ is essentially bounded from above and below if $\hat\psi$ is continuous and
the following two conditions hold:
\begin{enumerate}[(i)]
 \item $\exists\ \xi_0\neq 0,$ s.t. $\inf_{a\in[1,2]}|\hat\psi(a\xi_0)|>0$ 
 \item $\exists\ C>0$, s.t. $|\hat\psi(\xi)|^2\leq\frac{C|\xi|}{(1+|\xi|)^2},\ \forall \xi\in\mathbb{R}$
\end{enumerate}
 The canonical dual frame is another 
semi-discrete wavelet system. This can be seen if we use Corollary \ref{cor-trans-inv} and the observation that $m_{\Psi}(2^j\xi)=
m_{\Psi}(\xi)$ for all $j\in\mathbb{Z}$.
\begin{align*}
 S_\Psi^{-1} D_{2^j}\psi&=\mathcal{F}^{-1}\big(m_{\Psi}^{-1}D_{2^{-j}}\hat\psi\big)=\mathcal{F}^{-1}\big(m_{\Psi}^{-1}(2^j\ \cdot\ )
 D_{2^{-j}}\hat\psi\big)\\
 &=\mathcal{F}^{-1}\big(D_{2^{-j}}(m_{\Psi}^{-1}\hat\psi)\big) =D_{2^j}S_\Psi^{-1}\psi
 =:D_{2^j}\widetilde\psi
\end{align*}

 \begin{remark} \rm
Many common continuous transforms  from signal processing can be written as a translation invariant system
 and can therefore be treated with the previous results. Besides the STFT and the CWT,
the continuous shearlet transform \cite{kula12} and the continuous curvelet transform \cite{cado05} are noteworthy here.

It is important to observe that, although Theorem \ref{cnsgt} provides necessary and sufficient condition for a system of vectors 
on discrete LCA groups such as
$\mathbb{Z}$ or $\mathbb{C}^N$, to form 
a frame we did not investigate conditions for discrete frames on non-discrete LCA groups like $\mathbb{R}$ or $\mathbb{T}$.
\end{remark}
 
 \section{Reproducing pairs for the affine Weyl-\\ Heisenberg group}\label{aWH-section}

We will now further reduce the level of abstractness.
The affine Weyl-Heisenberg group $G_{aWH}$ and its unitary representations, see 
\cite{daforastte08,hola95,kato93,torr1,torr2}, is of particular relevance since it contains both, the 
affine group and the Weyl-Heisenberg group as subgroups which are the underlying groups of the continuous 
wavelet transform and the short-time Fourier transform. Furthermore, there is a wide range of transforms arising from this group
 and its subsets. Recently this transform has been successfully applied to medical data analysis \cite{dahtes08}.
 
Topologically,
$G_{aWH}$ is isomorphic to $\mathbb{R}^{2d}\times\mathbb{R}^\ast\times\mathbb{T}$ with the group law given by
\begin{center}
 $(x,\omega,a,\tau)\cdot(x',\omega',a',\tau')=(x+ax',\omega+\omega'/a,aa',\tau\cdot\tau'\cdot e^{-2\pi i\omega'\cdot x/a})$
\end{center}
and the neutral element  $e=(0,0,1,1)$.
The affine Weyl-Heisenberg group is unimodular and its Haar measure is given by $d\mu(x,\omega,a,\tau)=dxd\omega|a|^{-1}dad\tau$.
A unitary representation of $G_{aWH}$ on $L^2(\mathbb{R}^d)$ is given by
\begin{center}
 $\pi(x,\omega,a,\tau)\psi=\tau M_\omega T_x D_a \psi$
\end{center}
where the basic time-frequency operators on $\mathbb{R}^d$ are given by
\begin{center}
$T_x f(t)=f(t-x),$\ \ $M_\omega f(t)=e^{2\pi i\omega t}f(t)$,\ \  $D_a f(t)=|a|^{-d/2}f(t/a)$
\end{center}
Since $G_{aWH}$ is a locally compact group, the first step in the analysis of this representation is to examine if it is square-integrable.
Unfortunately, $\pi$ is not square-integrable because, loosely speaking, the group $G_{aWH}$ is too big.

To overcome this obstacle Torr{\'e}sani \cite{torr1} suggested to regularize the Haar measure by multiplying it with a weight function 
$\rho(\omega)$ and showed that under certain conditions this also leads to tight continuous frames. A different approach in the same paper 
considered subgroups of the affine Weyl-Heisenberg group to obtain square-integrability. For example, if $d=1$, the section
$(x,\eta_\lambda(a),a,\tau)$ with $\eta_\lambda(a)=\lambda\left(\frac{1}{a}-1\right)$ and $\lambda\in\mathbb{R}$ forms a
subgroup of $G_{aWH}$. Its representation is square-integrable
with left Haar measure $\frac{dxda}{|a|^2}$ and 
\begin{equation*}
 m_{\Psi,\Phi}\equiv\langle L\psi,L\varphi\rangle=\int_\mathbb{R}\overline{\hat\psi(\xi)}\hat\varphi(\xi)\frac{d\xi}{|\xi+\lambda|}
\end{equation*}
Within the scope of this paper we do not restrict ourselves to subgroups of $G_{aWH}$ but will apply the results from Section
\ref{cont-nonstationary gabor}. The key to reproducing pairs or continuous frames lies in an appropriate restriction of the group 
parameters and the choice of a measure $\mu$ on those subsets. Let
$\beta:\mathbb{R}^n\rightarrow \mathbb{R}$ and  $\eta:\mathbb{R}^n\rightarrow \mathbb{R}^d$, with $\beta,\eta$ piecewise continuous, 
$1\leq n\leq d$, and $\beta(\omega)=0$ only on a  null-set of $\mathbb{R}^n$. We
consider 
\begin{equation*}
G_{\beta,\eta}:=\Big\{(x,\eta(\omega),\beta(\omega),1): (x,\omega)\in \mathbb{R}^{d+n}\Big\}\subset G_{aWH}
\end{equation*}
together with the mapping
\begin{equation*}
\Psi(x,\omega):=M_{\eta(\omega)}T_x D_{\beta(\omega)}\psi,\ \ \psi\in L^2(\mathbb{R}^d)
\end{equation*}
$\Psi$ can be rewritten as 
\begin{equation*}\Psi(x,\omega)=e^{2\pi i x\cdot\eta(\omega)}\widetilde\Psi(x,\omega),\ \ 
\text{where}\ \ \widetilde\Psi(x,\omega):=T_x M_{\eta(\omega)}D_{\beta(\omega)}\psi
\end{equation*}
The new system $\widetilde \Psi$ is a nonstationary Gabor system. Hence, the recipe from Theorem 
\ref{cnsgt} is applicable to $\Psi$ by Lemma \ref{reprod-equiv-lemma}.

 The measure on $G_{\beta,\eta}$ will be defined by $d\mu_s(x,\omega):=|\beta(\omega)|^{s-d}dxd\omega$, $s\in\mathbb{R}$.
This particular choice of  
$d\mu_s$ is justified by two reasons. Firstly, the behavior of the system $\Psi$ is mainly depending on the scaling function $\beta$.
Secondly, the choice of $\beta$ and $\eta$ excludes certain choices of 
 $s$ if one wants to guarantee the existence of continuous frames. 

 To see this, let us assume that $d=1$ and consider the 
continuous wavelet transform, i.e., $\eta\equiv 0$, $\beta(\omega)=\omega$. In this case, one gets the symbol
$m_\Psi(\xi)=|\xi|^{-1-s}\int_\mathbb{R}|\hat\psi(a)|^2|a|^{s}da$. 
There is obviously no window $\psi\in L^2(\mathbb{R})$ such that this system forms a continuous frame if $s\neq-1$.
On the other hand, for the setup $\eta(\omega)=\omega$, 
$\beta(\omega)=(1+|\omega|)^{-1}$ no 
continuous frame exists if $s\neq 1$ which will be explained in more detail in Example 1.
These two short examples also indicate that, in most cases, there is no freedom in the choice of the parameter $s$ to obtain frames.

Theorem \ref{cnsgt} gives the following sufficient conditions for $(\Psi,\Phi)$ to form a reproducing system for $L^2(\mathbb{R}^d)$
\begin{equation}
 A\leq |m_{\Psi,\Phi}(\xi)|\leq B,\ \text{for\ a.e.}\ \xi\in\mathbb{R}^d
\end{equation}
where
\begin{equation}
m_{\Psi,\Phi}(\xi)=\int_{\mathbb{R}^n}\overline{\hat\psi\big(\beta(\omega)(\xi-\eta(\omega))\big)}\hat\varphi\big(\beta(\omega)
(\xi-\eta(\omega))\big)|\beta(\omega)|^{s}d\omega
\end{equation}
and 
\begin{equation}\label{abs-conv}
\int_{\mathbb{R}^n}\big|\hat\psi\big(\beta(\omega)(\xi-\eta(\omega))\big)\hat\varphi\big(\beta(\omega)(\xi-\eta(\omega))\big)\big|
 |\beta(\omega)|^{s}d\omega<\infty
\end{equation}
This result can be found in the literature for  particular choices of $\beta,\eta$.
In \cite{hola95} composite frames are investigated whereas the focus in \cite{daforastte08,dalo06} is on the
$\alpha$-transform and its uncertainty principles.
By choosing $\eta(\omega)=\omega$ we get a transform whose time-frequency resolution is frequency dependent. 
This is of particular interest for example in audio processing where one aims at constructing transforms following the time-frequency 
resolution of the human auditory system, see \cite{neccxxl13}.
\subsection{Example}\label{example-gawh-section}
 The following example has been introduced in \cite{hola95} where the notion of composite frames has been used.
 Let $d=n=1$, $\eta(\omega)=\omega$  and $\beta(\omega)=(1+|\omega|)^{-1}$. Substituting
 $z=\beta(\omega)(\xi-\omega)$ yields
 \begin{equation*}
 m_{\Psi,\Phi}(\xi)=\int_\mathbb{R}\overline{\hat\psi\Big(\frac{\xi-\omega}{1+|\omega|}\Big)}\hat\varphi\Big(\frac{\xi-\omega}{1+|\omega|}
 \Big)\frac{d\omega}{(1+|\omega|)^s}
  \end{equation*}
  \begin{equation*}
  =|1+\xi|^{1-s}\int_{-1}^\xi\overline{\hat\psi(z)}\hat\varphi(z)\frac{dz}{(1+z)^{2-s}}+
  |1-\xi|^{1-s}\int_{\xi}^1\overline{\hat\psi(z)}\hat\varphi(z)\frac{dz}{(1-z)^{2-s}}
 \end{equation*}
 It is not difficult to see that $ m_{\Psi,\Phi}$ fails to have either upper or lower bound if $s\neq1$. Hence we set $s=1$ and $ m_{\Psi,\Phi}$ 
 reads
  \begin{equation}\label{m-alpha-1}
 m_{\Psi,\Phi}(\xi)=\int_{-1}^\xi\overline{\hat\psi(z)}\hat\varphi(z)\frac{dz}{1+z}+\int_{\xi}^1\overline{\hat\psi(z)}\hat\varphi(z)
 \frac{dz}{1-z}
\end{equation}
This expression allows for explicit calculation of $m_{\Psi,\Phi}$ for many choices of windows $\psi,\varphi$. Take
 for instance $\hat\psi(\xi)=(1-\xi)\chi_{A}(\xi)$ and $\hat\varphi(\xi)=(1+\xi)\chi_{A}(\xi)$, with $A:=[-1,1]$, then 
 \begin{equation*}
  m_{\Psi,\Phi}(\xi)=\left\{\begin{array}{cl}2,& \mbox{for } |\xi|>1\\3-\xi^2, &\mbox{for } |\xi|\leq 1\end{array}\right.
 \end{equation*}
 In the following we will give an answer to question (iii) from Section \ref{prel} with the aid of the current example.

 \begin{proposition}
Let $\beta(\omega)=(1+|\omega|)^{-1},\ \eta(\omega)=\omega$ and $s=1$. There is no reproducing pair $(\Psi,\Phi)$, such that
\eqref{abs-conv} is satisfied and $S_{\Psi,\Phi}=Id$.
 \end{proposition}
\textbf{Proof:} Since \eqref{abs-conv} is satisfied,  Corollary \ref{cor-trans-inv} yields that $S_{\Psi,\Phi}=Id$
if and only if $m_{\Psi,\Phi}(\xi)=1$ for a.e. $\xi\in\mathbb{R}$. W.l.o.g. we may assume that $\hat\psi$ and $\hat\varphi$ are real-valued 
functions. Otherwise use $Re(\hat\psi\cdot\hat\varphi)$ instead of $\hat\psi\cdot\hat\varphi$ in the following arguments. Let us assume
that there exist  $\psi,\varphi\in L^2(\mathbb{R})$ s.t. $m_{\Psi,\Phi}=1$ for a.e. $\xi\in\mathbb{R}$. Then 
$m_{\Psi,\Phi}(\xi)=1$ for all $\xi\in[-1,1]$, since both summands in (\ref{m-alpha-1}) are continuous and consequently
$m_{\Psi,\Phi}'(\xi)=0$ for every  $\xi\in(-1,1)$. On the other hand, Lebesgue's differentiation theorem states that, for a.e. 
$\xi\in(-1,1)$, the derivative of  $m_{\Psi,\Phi}$ is given by
  \begin{equation*}
  m_{\Psi,\Phi}'(\xi)=\hat\psi(\xi)\hat\varphi(\xi)\left[\frac{1}{1+\xi}-\frac{1}{1-\xi}\right]=\frac{2\xi}{\xi^2-1}
  \hat\psi(\xi)\hat\varphi(\xi)
  \end{equation*}
Thus,
\begin{equation*}
\hat\psi(\xi)\hat\varphi(\xi)=0,\ \text{for a.e.}\ \xi\in(-1,1)
\end{equation*}
and $m_{\Psi,\Phi}(\xi)=0,\ \forall\xi\in[-1,1]$ which contradicts the assumption $m_{\Psi,\Phi}=1$, for a.e. $\xi\in\mathbb{R}$.
\hfill $\Box$\\ \\
Finally, we show that if $(\Psi,\Phi)$ is a reproducing pair, neither $\Psi$ nor $\Phi$ needs to be Bessel.
 To this end, take again $\hat\psi(\xi)=(1+\xi)\chi_{A}(\xi)$, $\hat\phi(\xi)=(1-\xi)\chi_{A}(\xi)$ and 
 $f \in L^1(\mathbb{R})\cap L^2(\mathbb{R})$ such that $|\hat f(\xi)|\geq 1,\ \forall \xi\in A$.
For every $(\xi,\omega)\in A\times\mathbb{R}$ it holds $\beta(\omega)(\xi-\omega)\in A$.
Hence, it follows
\begin{align*}
 \|C_\Psi f\|_{L^2(\mathbb{R}^2)}^2  &= \int_\mathbb{R}\int_\mathbb{R}|\hat f(\xi)|^2|\hat\psi(\beta(\omega)(\xi-\omega))|^2
\beta(\omega)d\xi d\omega\\
&\geq\int_\mathbb{R}\int_A|\hat\psi(\beta(\omega)(\xi-\omega))|^2\beta(\omega)d\xi d\omega\\
&=\int_\mathbb{R}\int_A\Big(1+\beta(\omega)(\xi-\omega))\Big)^2\beta(\omega)d\xi d\omega\\
&=\int_\mathbb{R}\int_A\big(1+2|\omega|\chi_{(-\infty,0]}(\omega)+\xi\big)^2\beta(\omega)^3d\xi d\omega\\
&=\int_\mathbb{R}\Big[C_1+\Big(C_2+2|\omega|\chi_{(-\infty,0]}(\omega)\Big)^2\ \Big]
\beta(\omega)^3 d\omega\\ &=\infty
\end{align*}
The same argument applies to $C_\Phi f$.
\section{Outlook and discussions}
It seems interesting to study discretization schemes for the CNSGT when starting from a semi-discrete system, see Section \ref{examples}.
Clearly, generalized coorbit theory \cite{fora05}
could be applied in this context. This approach however neither exploits that $Y$ is already discrete, nor
the abelian group structure of $G$.

Furthermore, as mentioned after Proposition \ref{sec:repordkern1}, a  characterization of the orbit of $C_\Psi$ is desired. 
To this end, a promising ansatz is to construct and investigate Gelfand triples of those spaces similar to the approach in 
\cite{jpaxxl09}. Moreover, as the images of $C_\Psi$ and $C_\Phi$ are mutually dual, an investigation in the context of partial 
inner product spaces \cite{antoine2009partial} seems worthwhile.

\section*{Acknowledgement}
This work was partially supported by the Austrian Science Fund (FWF) through the START-project 
FLAME (Frames and Linear Operators for Acoustical Modeling and Parameter Estimation): Y 551-N13.
They thank the anonymous reviewers for their input and Diana T. Stoeva for related discussions.

\bibliographystyle{abbrv}
\bibliography{paperbib}
\end{document}